\documentclass[12pt]{amsart}
\input{epsf.tex}
\usepackage{latexsym,amsfonts,amssymb,epsfig,verbatim}
\usepackage{amsmath, latexsym, graphics, textcomp}

\newcommand{\nc}{\newcommand}
\nc{\dmo}{\DeclareMathOperator}
\nc{\nt}{\newtheorem}

\theoremstyle{plain}
\nt{thm}{Theorem}[section]
\nt*{main}{Main Theorem}
\nt{lem}[thm]{Lemma}
\nt{cor}{Corollary}
\nt{prop}{Proposition}
\nt{fact}{Fact}

\def\Mod{\mbox{\rm{Mod}}}
\def\PMod{\mbox{\rm{PMod}}}
\def\Homeo{\mbox{\rm{Homeo}}}
\def\Aut{\mbox{\rm{Aut}}}
\def\Out{\mbox{\rm{Out}}}
\def\Hom{\mbox{\rm{Hom}}}

\nc{\bn}{B_n}
\nc{\pbn}{PB_n}
\nc{\cn}{A(B_{n-1})}
\nc{\hGa}{\hat{\Ga}}
\nc{\Q}{\mathbb{Q}}
\nc{\G}{\mathcal{G}}
\nc{\K}{\mathcal{K}}
\nc{\hG}{\hat{\G}}
\nc{\h}{\mathcal{H}}
\dmo{\fh}{H^1}
\nc{\invlim}{\ensuremath{\displaystyle{\lim_{\longrightarrow}\fh(\hG_i;\bz)}}}
\nc{\ilim}{\ensuremath{\displaystyle{\lim \fh(\hGa_i,\bz)}}}

\dmo{\aut}{Aut}
\dmo{\comm}{Comm}
\nc{\bz}{\mathbb{Z}}

\dmo{\tv}{Tv}
\dmo{\atv}{tv}

\nc{\modpms}{\Mod^\pm(S)}
\nc{\mods}{\Mod(S)}
\nc{\pmods}{\PMod(S)}

\nc{\dn}{D_n}

\nc{\bpf}{\begin{proof}} \nc{\epf}{\end{proof}}

\nc{\p}[1]{{\bf #1} }

\nc{\set}[1]{\{#1\}} \nc{\genby}[1]{\langle#1\rangle}

\nc{\beqn}{\begin{equation*}} \nc{\eeqn}{\end{equation*}}

\nc{\si}{\sigma} \nc{\De}{\Delta} \nc{\Ga}{\Gamma}

\begin{document}

\title{Abstract commensurators of braid groups}

\author{Christopher J Leininger}
\author{Dan Margalit}

\address{Department of Mathematics \\ Columbia University \\ 2990 Broadway MC 4448 \\
New York, NY 10027-6902}

\address{Department of Mathematics\\ University of Utah\\ 155 S 1440 East \\ Salt Lake City, UT 84112-0090}

\thanks{Both authors are supported by NSF postdoctoral fellowships.}

\email{clein@math.columbia.edu, margalit@math.utah.edu}



\maketitle

\begin{center}\today\end{center}

\begin{abstract}
Let $\bn$ be the braid group on $n \geq 4$ strands.  We show that
the abstract commensurator of $\bn$ is isomorphic to $\mods
\ltimes (\Q^\times \ltimes \Q^\infty)$, where $\mods$ is the
extended mapping class group of the sphere with $n+1$ punctures.
\end{abstract}

 \section{Introduction} \label{introsect}

Artin's \emph{braid group on $n$ strands}, denoted $\bn$, is the
group defined by the following presentation:
\begin{eqnarray*} \genby{ \si_1, \dots, \si_{n-1} : & \si_i \si_{i+1} \si_{i} = \si_{i+1} \si_i \si_{i+1} & \text{ for all } i, \\
& \si_i \si_j = \si_j \si_i & \text{ for } |i-j| > 1 }
.\end{eqnarray*}

This group also has a topological interpretation, from which it gets its name (see e.g. \cite{bb}).  We denote the center of $\bn$ by $Z$; it is infinite cyclic, generated by $z=(\si_1 \cdots \si_{n-1})^n$.  The goal of this paper is to characterize all isomorphisms between finite index subgroups of $\bn$.

The \emph{abstract commensurator} $\comm(G)$ of a group $G$ is the group of equivalence classes of isomorphisms of finite index
subgroups of $G$:
\[ \comm(G) = \set{\Phi: \Ga \stackrel{\cong}{\to} \De : \Ga,\De \text{ finite index subgroups of } G }/\sim \]
where $\Phi \sim \Phi'$  if there is a finite index subgroup
$\Ga'$ of $G$ such that $\Phi|_{\Ga'}=\Phi'|_{\Ga'}$.  The product of elements of $\comm(G)$ represented by $\Phi_1:\Ga_1\to\De_1$ and $\Phi_2:\Ga_2\to\De_2$ is an element represented by the isomorphism $\Phi_2 \circ \Phi_1|_{\Phi_1^{-1}(\De_1 \cap \Ga_2)}$.  A simple example is $\comm(\bz^n) \cong \textrm{GL}_n(\Q)$.

The \emph{extended mapping class group} $\mods$ of a surface $S$ is the group of isotopy classes of homeomorphisms of $S$:
\[ \mods = \pi_0(\Homeo^\pm(S)) .\]

\begin{main}Suppose $n \geq 4$, and let $S$ be the sphere with $n+1$ punctures.  Then we have:
\[ \comm(\bn) \cong \mods \ltimes \left( \Q^\times \ltimes \Q^\infty \right) .\]
\end{main}

We say that two groups $G$ and $G'$ are \emph{abstractly
commensurable} if they have isomorphic finite index subgroups.  In this case, it follows from the definition that $\comm(G) \cong
\comm(G')$.  A special case of this is when $G'$ is itself a
finite index subgroup of $G$.

Thus, the main theorem also gives the abstract commensurator of
all finite index subgroups of $\bn$.  We now mention two such
subgroups of general interest.  The \emph{pure braid group on n
strands} $\pbn$ is the kernel of the natural map from $\bn$ to the symmetric group on $\set{1,\dots,n}$ which sends $\sigma_i$ to the transposition switching $i$ and $i+1$.  The Artin group $\cn$ is isomorphic to the finite index subgroup of $\bn$ generated by $\si_1^2,\si_2,\dots,\si_{n-1}$.

\begin{cor} Suppose $n \geq 4$, and let $S$ be the sphere with $n+1$ punctures.  Then we have:
\[ \comm(\pbn) \cong \comm(\cn) \cong \mods \ltimes \left( \Q^\times \ltimes \Q^\infty \right) .\]
\end{cor}

The factor of $\mods$ in the main theorem comes from the following theorem of Charney and Crisp, which is a corollary of a theorem of Korkmaz \cite{cc} \cite{mk}.

\begin{thm} \label{kork}
Suppose $n \geq 4$, and let $S$ the sphere with $n+1$ punctures.
Then we have:
\[ \comm(\bn/Z) \cong \mods .\]
\end{thm}

This theorem relies on the classical fact that $\bn/Z$ is isomorphic
to the finite index subgroup of $\mods$ consisting of orientation
preserving elements that fix a single given puncture (see
e.g. \cite{cc}). Thus, there is a natural homomorphism $\mods \to \comm(\bn/Z)
\cong \comm(\mods)$, as $\mods$ acts on itself by inner
automorphisms.  Korkmaz's theorem is that this map is surjective, and
Charney and Crisp's contribution is that the map is injective (this is
implicit in the work of Ivanov \cite{nvi}).

\p{Other braid groups.} Our proof does not hold for $n=3$, as there is no analog of Theorem~\ref{kork}.  Indeed, $PB_3/Z$ is
isomorphic to the free group $F_2$, and $\comm(F_2)$ contains $\Aut(F_n)$ for all $n > 0$.  Also, note $B_2 \cong \bz$ and $B_1 = 1$.

\p{Historical background.} We think of $\comm(\bn)$ as describing ``hidden automorphisms'' of $\bn$ (compare with  \cite{nr} and \cite{fw}).
In that sense, our main result is a generalization of the theorem of Dyer and Grossman that $\Out(\bn) \cong \bz/2\bz$ \cite{dg}.  Recently, Charney and Crisp proved that $\Out(\cn) \cong (\bz \rtimes \bz/2\bz) \times \bz/2\bz$ \cite{cc}, and Bell and Margalit proved that $\Out(\pbn) \cong \bz^N \rtimes (\Sigma_n \times \bz/2\bz)$, where $\Sigma_n$ is the symmetric group on $n$ letters \cite{bm2}.  Charney and Crisp also showed that the abstract commensurator of any finite type Artin group (e.g. $\bn$) contains an infinitely generated free abelian subgroup.  Very recently, Crisp has determined the abstract commensurators of certain 2-dimensional Artin groups \cite{jc}.

As explained, there is a close connection between braid groups and mapping class groups.
There are several recent results giving the abstract commensurators of subgroups of $\mods$.
Ivanov showed that $\comm(\mods)$ is isomorphic to $\mods$ for most surfaces $S$ \cite{nvi}; the genus zero case is Theorem~\ref{kork}.  Farb and Ivanov proved that the abstract commensurator of the Torelli subgroup of $\mods$ is isomorphic to $\mods$ \cite{fi}.
Most recently, Brendle and Margalit showed that the so-called Johnson kernel, a subgroup of the Torelli group, also has abstract commensurator $\mods$ \cite{brm}.  A similar phenomenon exists with the related group $\Out(F_n)$; it is a result of Farb and Handel that $\comm(\Out(F_n))$ is isomorphic to $\Out(F_n)$ when $n \geq 4$ \cite{fh}.

Some of the ingredients in our proof can be viewed as generalizations of facts about automorphism groups (see Section \ref{generalsection}), and are likely well known to others familiar with commensurators, though we have found no references for them.

\p{Outline of proof.} First we find a group $\G$ which is abstractly
commensurable to $\bn$, and which is a direct product over its center, $Z$.

\emph{Proposition 1}: $\G = \hG \times Z$.

We then define the transvection subgroup $\tv(\G)$ of $\comm(\G)$
and show this group splits off as a semidirect factor.

\emph{Proposition 2}: $\comm(\G) \cong \mods \ltimes \tv(\G)$.

To understand the structure of $\tv(\G)$, we define the subgroup $\h$ of simple transvections and show this groups splits from $\tv(\G)$ as a semidirect factor.

\emph{Proposition 3}: $\tv(\G) \cong \Q^\times \ltimes \h$.

Finally, we use the notion of a divisible group to describe $\h$.

\emph{Proposition 4}: $\h \cong \Q^\infty$.

\p{Acknowledgements.} We are grateful to Bob Bell, Mladen Bestvina,
Joan Birman, Benson Farb, Walter Neumann, and Kevin Wortman for much
encouragement and many enjoyable conversations.  The second author
would like to thank the mathematics department of Columbia University
for providing a very pleasant and stimulating environment for the
visit during which this project was begun.


\section{The proof}

Let $n \geq 4$ be fixed.  We start by finding a group $\G$ which is
abstractly commensurable to $\bn$, and
which splits over its center.  A priori, this
is an easier group to work with than $\bn$, and as mentioned earlier,
it has an isomorphic abstract commensurator.

\p{Length homomorphism.} We will make use of the \emph{length homomorphism} $L:\bn \to
\bz$, which is defined by $\si_i \to 1$ for all $i$.  Note that
$L$ is indeed a homomorphism, and that $L(z)=n(n-1)$.

\begin{prop} \label{p1}
$\bn$ is abstractly commensurable to the external direct product \[ \G = \hG \times Z \]
where $\hG$ is a finite index subgroup of $\bn/Z$.\end{prop}

\bpf

Let $\K$ be the kernel of the composition
\[ \bn \stackrel{L}{\to} \bz \to \bz/n(n-1)\bz \]
where the latter map is reduction modulo $n(n-1)$.

Since $Z < \K$, we have:
\begin{equation}\label{seszggz}
1 \to Z \to \mathcal{K} \to \K/Z \to 1 .
\end{equation}
We can view the restriction of $L$ to $\K$ as a projection to $Z$ by
defining a map $\K \to Z$ via $g \mapsto z^{L(g)/n(n-1)}$; this is a splitting for the sequence.
Thus, $\K$ is isomorphic to the external direct product which we denote $\G = \hG \times Z \cong \K$ where $\hG = \K/Z$.\epf

As an abuse of notation, we will identify $\hG$ and $Z$ with their images in $\G$.



\p{Transvections.} We define the \emph{transvection subgroup}
$\tv(\G)$  of $\comm(\G)$ by the following short exact sequence (compare with
\cite{cc}):
\begin{equation} \label{sescomzggz}
1 \to \tv(\G) \to \comm(\G) \to \comm(\hG) \to 1
.\end{equation}
That $\comm(\G)$ surjects onto $\comm(\hG)$ follows directly from the fact that $\G = \hG \times Z$.

\begin{prop} $\comm(\G) \cong \mods \ltimes \tv(\G)$.\end{prop}

\bpf

Since $\hG$ is a finite index subgroup of $\bn/Z$,
Theorem~\ref{kork} gives $\comm(\hG) \cong \mods$.  We define a
splitting for (\ref{sescomzggz}) by sending an element of
$\comm(\hG)$ represented by $\Psi:\hGa \to \hat{\De}$ to the element
of $\comm(\G)$ represented by $\Psi \times
1: \hGa \times Z \to \hat{\De} \times Z$.\epf



\p{Virtual center of $\G$.} We will use the fact that the structure of $\G$ with respect to its center is preserved under passage to finite index subgroups:

\begin{lem}\label{virtcent}If $\Ga$ is any finite index subgroup of $\G$, then $Z(\Ga)=\Ga \cap Z$.\end{lem}

Via the isomorphism of Proposition~\ref{p1}, this lemma is part of the
proof of Theorem~\ref{kork}; it is equivalent to the statement that
the map from $\mods$ to $\comm(\mods)$ is injective.

\p{Simple transvections.} There is a homomorphism $\theta$ from $\tv(\G)$ to $\Q^\times$ which measures the action on $Z$.  Indeed, given an element of $\tv(\G)$ represented by $\Phi:\Ga \to \De$, and any $z^q \in \Ga$, we must have $\Phi(z^q)=z^p$ for some nonzero $p$ (by Lemma~\ref{virtcent}); define $\theta([\Phi])$ to
be $p/q$.  This is a well-defined homomorphism; we call its kernel the group $\h$ of \emph{simple transvections}:
\begin{equation} \label{sestvzggz}
1 \to \h \to \tv(\G) \stackrel{\theta}{\to} \Q^\times \to 1
.\end{equation}

\begin{prop} $\tv(\G) \cong \Q^\times \ltimes \h$ .\end{prop}

\bpf

There is a splitting of (\ref{sestvzggz}):
given $p/q \in \Q^\times$, where $p,q \in \bz$, let $\Phi: \hG \times \genby{ z^q } \to \hG \times \genby{ z^p }$ be the transvection which is the identity on the first factor, and sends $z^q$ to $z^p$.\epf



\p{Subgroup structure.}  Given finite index subgroups $\hGa < \hG$ and $Z_0 < Z$, by an abuse of notation, we identify the external direct product $\hGa \times Z_0$ with its image under the obvious inclusion $\hGa \times Z_0 < \hG \times Z$.
The next lemma says that we can always choose subgroups of this type as domains for representatives of commensurators.

\begin{lem} \label{ss}If $\Ga$ is a finite index subgroup of $\G$, then $\Ga$ has a finite index subgroup $\Ga'$ of the form:
\[ \Ga' = \hGa' \times Z(\Ga) \]
where $\hGa'$ is a finite index subgroup of $\hG$.\end{lem}

\bpf

By Lemma~\ref{virtcent}, $Z(\Ga) < Z$.
Let $\Ga'$ be the kernel of the composition
\[ \Ga \hookrightarrow \G \to Z \to Z/Z(\Ga) \]
where the latter two maps are the obvious projections.  Denote the
composition of the first two maps by $\pi$.  Then the short exact
sequence
\[ 1 \to Z(\Ga) \to \Ga' \to \Ga'/Z(\Ga) \to 1 \]
has a splitting $\Ga' \to Z(\Ga)$ given by $g \mapsto \pi(g)$ with kernel $\hGa' < \hG$, and the lemma follows. \epf


\p{Transvections and cohomology.} In order to get a clearer picture of $\h$, we need a description of elements of $\tv(\G)$.  Recall that for a group $G$, we have $\fh(G,\bz) \cong \Hom(G,\bz)$.

\begin{lem}\label{tvchar}
Suppose $\Phi:\Ga \to \De$ represents an element of $\tv(\G)$.
Then there exists $\phi \in \fh(\Ga,\bz)$ so that $\Phi$ is given by:
\[ \Phi(g) = gz^{\phi(g)} .\]
If $[\Phi] \in \h$ and $\Ga \cong \hGa \times Z(\Ga)$, we may view $\phi$ as an element of $\fh(\hGa,\bz)$.
\end{lem}

By Lemma~\ref{ss}, the second statement applies to all elements of
$\h$.

\bpf

We define an element $\phi \in \fh(\Ga,\bz)$ by the equation
\[ z^{\phi(g)} = g^{-1} \Phi(g) .\]
That $\phi$ is a homomorphism follows from the assumption that $\Phi$ is a transvection (in particular, $g^{-1} \Phi(g)$ is central):
\[ z^{\phi(gh)} = h^{-1} g^{-1} \Phi(g) \Phi(h) = g^{-1} \Phi(g) h^{-1} \Phi(h) = z^{\phi(g)+\phi(h)} \]
The first statement follows.  The second statement is clear:  if $[\Phi] \in \h$, then $\Phi(g) = g$; thus $\phi|_{Z(\Ga)}=0$, and $\phi$ descends to $\hGa \cong \Ga / Z(\Ga)$.\epf


\p{Direct limits.} Let $I$ be a directed partially ordered set; that
is, $I$ is a partially ordered set with the property that for any $i,j
\in I$, there is a $k \in I$ with $i,j \leq k$.  A collection
of abelian groups $\set{G_i}_{i \in I}$ and homomorphisms $\set{
f_{ij}: G_i \to G_j }_{\set{i,j \in I | i \leq j} }$ forms a
\emph{direct system} if: (1) the homomorphism $f_{ii}$ is the identity
map for all $i$; and (2) given any two
homomorphisms $f_{ij}$ and $f_{jk}$, we have $f_{ik} = f_{jk} \circ f_{ij}$.

The \emph{direct limit} of the direct system $(G_i,f_{ij})$, which we
denote $\lim G_i$, is the group which satisfies the following
universal property: if $\iota_{i}:G_{i} \rightarrow G$ is a collection
of homomorphisms respecting the homomorphisms $f_{ij}$, then there is
a unique homomorphism $\iota: \lim G_i \rightarrow G$ through which
each $\iota_i$ factors.  It follows that if each $\iota_i$ is a
monomorphism, then so is $\iota$.  In this case, each $G_{i}$
naturally includes into $\lim G_i$.

\p{Direct limit of cohomology groups.} We will consider the direct
system of groups $\fh(\hGa_i,\bz)$, where $\hGa_i$ ranges over all
finite index subgroups of $\hG$.  Since $\bz$ is torsion free, the
natural homomorphism $\fh(\hGa_1,\bz) \to \fh(\hGa_2,\bz)$, for any
$\hGa_2 < \hGa_1$, is injective; the inclusion is given by restriction
of homomorphisms.  That (the index set of) $\set{\hGa_i}$ forms a
directed partially ordered set is the fact that any two finite index
subgroups have a common finite index subgroup (namely, their
intersection).  Properties (1) and (2) of direct systems are apparent.
Thus, $\ilim$ is defined.

\begin{lem}\label{invlim} $\h \cong \ilim$.\end{lem}

\bpf

If $\hGa_{i}$ is any finite index subgroup of $\hG$, then there is a monomorphism
\[ \Xi_{\hGa_{i}} : \fh(\hGa_{i},\bz) \to \h \]
given by
\[ \phi \mapsto \left[ \Phi:\hGa_{i} \times Z \stackrel{\cong}{\to} \Ga' \right] \]
where $\Phi(g) = gz^{\phi(g)}$.

Moreover, this map respects the inclusions of the aforementioned direct system: if $\hGa_2 < \hGa_1$ then $\fh(\hGa_1,\bz) < \fh(\hGa_2,\bz)$ by restriction.
Thus $\Xi_{\hGa_1}$ is the restriction of $\Xi_{\hGa_2}$, and so the universal property guarantees a well-defined  injection $\Xi: \ilim \hookrightarrow \h$.
The inverse of $\Xi$ is provided by Lemma~\ref{tvchar}, so $\Xi$ is an isomorphism.\epf


\p{Divisible groups.} An abelian group $G$ is \emph{divisible} if for any element $g$ of $G$, and any positive integer $q$, there is an element $h$ of $G$ with $h^q=g$.  The next fact follows from definitions:

\begin{fact}\label{stdg} Any torsion free divisible group is a vector space over $\Q$.\end{fact}

\p{Braid groups and free groups.} Given the inclusion of $\bn/Z$ into
$\mods$ described in Section~\ref{introsect}, it follows from the
definition of $\pbn$ that $\pbn/Z$ is isomorphic to the subgroup of
$\mods$ consisting of orientation preserving mapping classes that fix
each puncture (note $Z(\pbn)=Z$).  Thus, for any $m < n$, there is a
surjection $PB_n/Z \to PB_m/Z(PB_m)$ obtained by ``forgetting'' $n-m$
of the punctures.  In this way, $\pbn/Z$ maps surjectively onto
$PB_3/Z(PB_3)$, which is isomorphic to the free group on two letters
$F_2$ (the last statement follows, for example, from the Birman exact
sequence \cite[Theorem 1.4]{jb}).  We record this fact for future
reference:

\begin{fact}\label{freebraid} If $n \geq 3$, then $PB_n/Z$ surjects onto $F_2$.\end{fact}

The following proposition completes the proof of the main theorem.

\begin{prop}\label{qi} $\h \cong \Q^\infty$.\end{prop}

\bpf

By Lemma~\ref{invlim}, we need only prove $\ilim \cong \Q^\infty$.
By Fact~\ref{stdg}, and since $\ilim$ is countable, it suffices to show that $\ilim$ is a torsion free divisible group which contains free abelian subgroups of arbitrary rank.

First, $\ilim$ is torsion free and abelian since each $\fh(\hGa_i,\bz)$ has these properties.

To see that $\ilim$ is a divisible group, let $\phi \in \fh(\hGa_i,\bz) \subset \ilim$, and let $q \in \bz$.
Consider the subgroup $\hGa_q$ of $\hGa$ that is the kernel of the composition: \[ \hGa \stackrel{\phi}{\to} \bz \to \bz/q\bz .\] Then $\phi|_{\hGa_q}$ maps to $q\bz$.  Thus, $\phi = q \phi'$ for some $\phi' \in \fh(\hGa_q,\bz)$.

We now construct free abelian subgroups of arbitrarily large rank in $\ilim$.
Let $\Pi:\pbn/Z \to F_2$ be the surjection given by Fact~\ref{freebraid}.
We can choose finite index free subgroups $F_k < F_2$ with any rank $k \geq 2$.
Thus, since surjections and passages to finite index subgroups both induce inclusions on cohomology, we have the required injections:
\[ \bz^k \cong \fh(F_k,\bz) \hookrightarrow \fh(\Pi^{-1}(G_i) \cap
\hG,\bz) \hookrightarrow \ilim \]
\epf



\section{Generalities} \label{generalsection}

As we have mentioned, given a group $G$, we view $\comm(G)$ as a generalization of $\aut(G)$.
Automorphism groups of central extensions of centerless groups can be understood as follows.
First, let $G$ be a group with $Z(G) = 1$, $A$ an abelian group, and
\begin{equation} \label{sesgc1}
1 \to A \to \Ga \to G \to 1
\end{equation}
a split central extension.
This induces a split exact sequence
\begin{equation} \label{sesgc2}
1 \to \atv(\Ga) \to \aut(\Ga) \to \aut(G) \to 1.
\end{equation}
The subgroup  $\atv(\Ga)$ of $\aut(\Ga)$ consists of those automorphisms which become trivial upon passing to the quotient $G$.
This group fits into a split exact sequence
\begin{equation} \label{sesgc3} 1 \to \fh(G,A) \to \atv(\Ga) \to \aut(A) \to 1.
\end{equation}
The inclusion of $\fh(G,A)$ into $\atv(\Ga)$ is defined by
sending any $\phi$ in $\fh(G,A)$ to the map given by $g \mapsto g\tilde{\phi}(g)$, where $\tilde{\phi}$ is the pullback of $\phi$ to $\fh(\Ga,A)$ (compare Lemma~\ref{tvchar}).

If (\ref{sesgc1}) is not split, then we still have sequences (\ref{sesgc2}) and (\ref{sesgc3}), but these need not be exact (although the second map in each is still injective).

When all finite index subgroups of $G$ are centerless, we obtain a completely analogous picture for $\comm(G)$.
In particular, (\ref{sesgc2}) becomes
\begin{equation*}\tag{\ref{sesgc2}'} \label{sesgc2'}
1 \to \tv(\Ga) \to \comm(\Ga) \to \comm(G) \to 1.
\end{equation*}
The group $\tv(\Ga)$ consists of those commensurators which are trivial in $\comm(G)$.
This also determines a sequence analogous to (\ref{sesgc3}):
\begin{equation*}\tag{\ref{sesgc3}'} \label{sesgc3'}
1 \to \lim\fh(G_{i},A) \to \tv(\Ga) \to \comm(A) \to 1.
\end{equation*}
When (\ref{sesgc1}) virtually splits, these are also split exact.
Otherwise, they need not be exact.

Finally, we have $\lim \fh(G_{i},A) \cong \lim \fh(G_i,A/T)$, where
$T$ is the torsion subgroup of $A$.  If $A$ is finitely generated and
$A/T \cong {\mathbb Z}^{m}$, then
$$\lim \fh(G_{i},A) \cong \Q^{m \cdot vb_{1}} \mbox{ and } \comm(A) \cong \textrm{GL}_m(\Q)$$
where $vb_{1} \in {\mathbb Z}_{\geq 0} \cup \{ \infty \}$ is the {\em
  virtual first betti number of $G$}, i.e. the supremum of the first
betti numbers of finite  index subgroups of $G$.  So in this case, we
have:
\[ \comm(\Ga) \cong \comm(G) \ltimes (\textrm{GL}_m(\Q) \ltimes \Q^{m
  \cdot vb_1}) \]


\bibliographystyle{plain}
\bibliography{commbn}

\end{document}